\documentclass[twoside,leqno]{article}

\usepackage[letterpaper]{geometry}

\usepackage{siamproceedings}
\usepackage{mathtools}
\usepackage{amsmath}
\usepackage{amssymb}
\usepackage{stmaryrd}

\usepackage[T1]{fontenc}
\usepackage{amsfonts}
\usepackage{graphicx}
\usepackage{epstopdf}
\usepackage{enumitem}
\usepackage{algorithmic}
\ifpdf
  \DeclareGraphicsExtensions{.eps,.pdf,.png,.jpg}
\else
  \DeclareGraphicsExtensions{.eps}
\fi

\usepackage{tikz-cd}

\usetikzlibrary{decorations.markings}
\usetikzlibrary{decorations.pathmorphing}
\tikzset{
  module/.style={
    postaction={decorate},
    decoration={
      markings,
      mark=at position #1 with {\arrow{|}}}},
  module/.default=0.5,
  we/.style=
  { postaction={%
      decorate,
      decoration={
        markings,
        mark=at position #1 with {%
          \node[transform shape, yshift=.2em]{%
            \resizebox{0.5em}{!}{$\sim$}};}}}},
  we/.default=0.5,
  we'/.style=
  { postaction={%
      decorate,
      decoration={
        markings,
        mark=at position #1 with {%
          \node[transform shape, yshift=-.2em, rotate=180]{%
            \resizebox{0.5em}{!}{$\sim$}};}}}},
  we'/.default=0.5,
  iso/.style=
  { postaction={%
      decorate,
      decoration={
        markings,
        mark=at position #1 with {%
          \node[transform shape, yshift=.2em]{%
            \resizebox{0.5em}{!}{$\simeq$}};}}}},
  iso/.default=0.5,
  iso'/.style=
  { postaction={
      decorate,
      decoration={
        markings,
        mark=at position #1 with {%
          \node[transform shape, yshift=-.2em, rotate=180]{%
            \resizebox{0.5em}{!}{$\simeq$}};}}}},
  iso'/.default=0.5,
}

\tikzset{
  proarrow/.style={->, module},
  proequal/.style={-, double, module},
  prodotted/.style={->,dotted, module},
  prodashed/.style={->,dashed, module},
  wearrow/.style={->, we},
  wedashed/.style={->, dashed, we},
  wedotted/.style={->, dotted, we},
  tfibarrow/.style={->>, we=0.45},
  tfibdotted/.style={->>,dotted, we=0.45},
  tfibdashed/.style={->>,dashed, we=0.45},
  tcofarrow/.style={>->, we},
  tcofdashed/.style={>->, dashed, we},
  uwearrow/.style={->, we'},
  uwedashed/.style={->, dashed, we'},
  uwedotted/.style={->, dotted, we'},
  utfibarrow/.style={->>, we'=0.45},
  utfibdotted/.style={->>,dotted, we'=0.45},
  utfibdashed/.style={->>,dashed, we'=0.45},
  utcofarrow/.style={>->, we'},
  isoarrow/.style={->, iso},
  isodashed/.style={->, dashed, iso},
  uisoarrow/.style={->, iso'},
  uisodashed/.style={->, dashed, iso'},
  isocell/.style={=>, iso},
  isocelldashed/.style={=>, dashed, iso},
  uisocell/.style={=>, iso'},
  uisocelldashed/.style={=>, dashed, iso'}
}

\makeatletter
\def\makeslashed#1#2#3#4#5{#1{\mathpalette{\sla@{#2}{#3}{#4}}{#5}}}

\def\@mathlower#1#2#3{\setbox0=\hbox{$\m@th#2#3$}\lower#1\ht0\box0}
\def\mathlower#1#2{\mathpalette{\@mathlower{#1}}{#2}}
\makeatother

\newcommand\tailxrightarrow[2][]{%
  \mathrel{\ooalign{$\xrightarrow[#1\mkern4mu]{#2\mkern4mu}$\cr%
  \hidewidth$\Yright\mkern17mu$}}
}

\newcommand{\fto}{\twoheadrightarrow}

\newcommand{\cwto}{\tailxrightarrow{{\smash{\mathlower{0.8}{\sim}}}}}

\newsiamremark{remark}{Remark}
\newsiamremark{construction}{Construction}
\crefname{construction}{Construction}{Construction}
\newsiamthm{claim}{Claim}
\newsiamremark{notation}{Notation}
\crefname{notation}{Notation}{Notation}
\newsiamremark{example}{Example}
\crefname{example}{Example}{Example}
\newsiamremark{warning}{Warning}
\crefname{warning}{Warning}{Warning}

\newcommand{\comp}{\textup{comp}}
\newcommand{\colim}{\textup{colim}}
\newcommand{\ev}{\textup{ev}}
\newcommand{\evid}{\textup{ev}_{\textup{id}}}
\newcommand{\Fun}{\textup{Fun}}
\newcommand{\Hom}{\textup{Hom}}
\newcommand{\id}{\textup{id}}
\newcommand{\Map}{\textup{Map}}
\newcommand{\op}{\textup{op}}
\newcommand{\arrowind}{\textup{arrow-ind}}
\newcommand{\pathind}{\textup{path-ind}}
\newcommand{\patheq}{\textup{path-to-eq}}
\newcommand{\arrfun}{\textup{arr-to-fun}}
\newcommand{\refl}{\textup{refl}}
\newcommand{\tr}{\textup{tr}}

\newcommand{\type}[1]{{\mathsf{#1}}}
\newcommand{\univ}{\mathcal{U}}
\newcommand{\cS}{\mathcal{S}}
\newcommand{\Aut}{\type{Aut}}
\newcommand{\Fin}{\type{Fin}}

\newcommand{\Eq}{\mathsf{Eq}}
\newcommand{\isEquiv}{\mathsf{isEquiv}}
\newcommand{\isIso}{\mathsf{isIso}}
\newcommand{\isContr}{\mathsf{isContr}}
\newcommand{\fib}{\mathsf{fib}}

\DeclareMathAlphabet{\mathbbe}{U}{bbold}{m}{n}

\newcommand{\1}{\mathbbe{1}}
\newcommand{\2}{\mathbbe{2}}
\newcommand{\3}{\mathbbe{3}}
\newcommand{\iso}{\mathbb{I}}
\newcommand{\DDelta}{\mathbbe{\Delta}}

\newcommand{\cat}[1]{\textup{\textsf{#1}}}
\newcommand{\fun}[1]{\textup{#1}}

\newcommand{\Set}{\mathcal{S}\cat{et}}
\newcommand{\sSet}{\cat{s}\mathcal{S}\cat{et}}
\newcommand{\Top}{\mathcal{T}\cat{op}}

\newcommand{\coslice}[2]{{{}^{{#1}/}\!{#2}}}

\makeatletter
\let\c@equation\c@theorem
\makeatother
\numberwithin{equation}{section}

\font\maljapanese=dmjhira at 2ex 
\def\yo{\textrm{\maljapanese\char"48}}

\usepackage{amsopn}

\begin{document}

\newcommand\relatedversion{}

\title{\Large Synthetic perspectives on spaces and categories\relatedversion}
    \author{Emily Riehl\thanks{This article is dedicated to Jack Morava, whose generous spirit and playful intellect are sorely missed (\email{eriehl@jhu.edu}, \url{emilyriehl.github.io}).}}

\date{}

\maketitle


\fancyfoot[R]{\scriptsize{Copyright \textcopyright\ 2026 by SIAM\\
Unauthorized reproduction of this article is prohibited}}





\begin{abstract} Recently discovered domain-specific formal systems---specifically \emph{homotopy type theory} and \emph{simplicial type theory}---provide new perspectives on spaces and categories in a natively equivalence-invariant setting. In this note, we expose fundamental proof techniques from these parallel settings: describing induction principles over paths or arrows and constructions involving universes that are either univalent or directed univalent. 
\end{abstract}

\section{Introduction.}

Over the past several decades, various formalisms have been developed for working with spaces or categories in a natively ``derived'' or equivalence-invariant manner \cite{Grothendieck1957, Verdier1996, Quillen1967, Freyd1976pi, Makkai1995fo, Joyal2002,  DwyerHirschhornKanSmith2004, Lurie2009}. 

\subsection{On spaces.}\label{ssec:on-spaces}

In the case of topological spaces, the notion of equivalence we have in mind is \emph{weak homotopy equivalence}, a map $f \colon X \to Y$ that induces an isomorphism on path components and all homotopy groups. Constructions that respect this notion of equivalence are constructions on \emph{homotopy types}---i.e., spaces up to weak homotopy equivalence---which are also known as $\infty$-\emph{groupoids} (following Grothendieck \cite{Grothendieck2022ps}) or \emph{anima} (following \v{C}esnavi\v{c}ius, Clausen, and Scholze \cite[\S 11.1]{ClausenScholze2020}, \cite[\S 1.2]{CesnaviciusScholze2024}). In this note, we use the more familiar term ``spaces'' to refer to these equivalent notions but caution that our spaces do not have a well-defined underlying set: instead, a point in a space refers to a homotopy class of maps out of a contractible space, i.e., an element of the set of path components. 

Our aim in the first part of this note is to provide a prospective on spaces that captures the insights of \emph{homotopy type theory} and \emph{univalent foundations}, a ``surprising synthesis of constructive intensional type theory and abstract homotopy theory'' \cite{Shulman2017cwp} discovered in the early 21st century. Homotopy type theory can be thought of as a synthetic theory of spaces, where spaces are primitive notions, usually called ``(homotopy) types.'' We will introduce notations, constructions, and proofs inspired by that formal framework in a concrete setting that can be described in traditional foundations. 

To make sense of the constructions that follow, we must work in a ``very convenient category of spaces.'' In \cite{Steenrod1967}, Steenrod explained how to construct a category that contained the familiar examples of spaces and had the property of being \emph{cartesian closed}, meaning that the cartesian product with a fixed space $X$ has a right adjoint, defining an internal hom: for spaces $X$ and $Y$, this constructs the \emph{mapping space} $\Map(X,Y) \cong Y^X$ from $X$ to $Y$. Here we require our category of spaces to be \emph{locally cartesian closed}, a notion we review in \S\ref{sec:spaces}. In fact our intended semantic setting is more general than our language will suggest. Shulman has shown \cite{Shulman2019} that homotopy type theory, and in particular all of our constructions, can be interpreted in an arbitrary \emph{$\infty$-topos} \cite{Rezk2005, Lurie2009}, of which the base example is the $\infty$-topos of spaces. Thus our primitive ``spaces'' might be somewhat more elaborate.

\subsection{On categories.}

Of the many interesting variants of the notion of a category, here we focus on $(\infty,1)$-\emph{categories}---first introduced by Boardman and Vogt \cite{BoardmanVogt1973}, further developed by Joyal \cite{Joyal2002}, and nicknamed ``$\infty$-\emph{categories}'' by Lurie \cite{Lurie2009}---which are the analogs of ordinary categories replacing all the sets in the definition by spaces, in the sense introduced above. This can be understood as \emph{enriching} the ordinary notion of a category---adding morphisms between morphisms, and morphisms between morphisms, and morphisms between morphisms between morphisms ad infinitum thought of as paths, homotopies, and higher homotopies in the spaces of arrows. But on account of the fact that all of our constructions on spaces are required to be equivalence-invariant, this must also be understood as \emph{weakening} the original notion. The central structure in an ordinary strict 1-category is given by its composition function, defined for any triple of objects $x$, $y$, $z$: \[ \begin{tikzcd} \Hom(y,z) \times \Hom(x,y) \arrow[r, "\circ"] & \Hom(x,z) \end{tikzcd}\]  
Here, this takes the form of a map between equivalence-invariant spaces of arrows, and consequently the composition of morphisms $g \colon Y \to Z$ and $f\colon X \to Y$ is only well-defined up to a contractible space of choices, which is the homotopical analog of the usual requirement that composition is defined uniquely. In this note, we use the term ``categories'' to refer to $\infty$-categories and use the term ``(strict) 1-categories'' to refer to the ordinary notion in the few settings it appears.

As categories, in our sense, are built out of spaces, it is attractive to develop the theory of categories in the synthetic setting for spaces described in \S\ref{ssec:on-spaces}. This avoids the notorious complexity inherent to the traditional ``analytic'' definitions of categories 
surveyed in \cite{AntolinCamarena2016, Bergner2018}.
To the formal system of homotopy type theory, we add an axiomatic directed arrow satisfying strict interval axioms. In \S\ref{sec:cats}, we sketch the structures present in our setting for synthetic category theory, known in the literature as \emph{simplicial type theory} \cite{RiehlShulman2017}. Again, our intended semantic setting is more general than our language suggests. All of our constructions will make sense not just for categories---the $\infty$-categories described above---but to what are called \emph{internal $\infty$-categories}, meaning $\infty$-categories defined in an arbitrary $\infty$-topos \cite{Martini2022yo,MartiniWolf2024}, or in even more general settings \cite{Rasekh2025sc}.

\subsection{Outline.}\label{ssec:plan} We could easily fill this article with more precise descriptions of the formal frameworks for synthetic spaces and synthetic categories---homotopy type theory and simplicial type theory respectively. Instead, our plan is to highlight some of the fundamental constructions that can equally be unpacked in traditional foundations. In each parallel setting, we explore
 new induction principles, specifically \emph{path induction} in \S\ref{sec:path} and \emph{arrow induction} in \S\ref{sec:arrow}, and
   the use of universes that are respectively \emph{univalent} in \S\ref{sec:ua} and \emph{directed univalent} in \S\ref{sec:dua}.
 In \S\ref{sec:vistas}, we conclude by highlighting ongoing work of our colleagues concerning synthetic spaces, synthetic categories, and synthetic mathematics more broadly.

There is a further benefit of synthetic mathematical frameworks that will not feature here but also motivates this work: namely their amenity to computer formalization, as formal proofs in synthetic domains tend to more closely parallel their pen-and-paper counterparts. In particular, theorems in homotopy type theory, such as the applications of univalence mentioned in \S\ref{ssec:ua-applications} have been formalized many times over, while the results on synthetic categories appearing in \S\ref{sec:cats}--\ref{sec:arrow} have been formalized in the computer proof assistant \textsc{Rzk} that implements the simplicial type theory \cite{Kudasov2023rzk,KudasovRiehlWeinberger2024, sHoTT}, the first and to date essentially only formalization of results from higher category theory. This proof assistant did not exist while the paper \cite{RiehlShulman2017} was being written, so the formalization was undertaken after the fact, catching the errors in the original paper noted in \S\ref{sec:arrow}.

\section{A very convenient category of spaces.}\label{sec:spaces}

Following Eilenberg and Zilber \cite{EilenbergZilber1950}, Gabriel and Zisman \cite{GabrielZisman1967}, Kan \cite{Kan1957css}, Quillen \cite{Quillen1967} and others, the classical homotopy theory of homotopy types can be modelled by a particular strict presheaf category, namely the strict 1-category $\sSet \coloneqq \Set^{\DDelta^\op}$ of \emph{simplicial sets}, indexed by the topologists' strict 1-category $\DDelta$ of finite non-empty ordinals and order-preserving maps. The representables---the objects in the image of the Yoneda embedding $\yo \colon \DDelta \hookrightarrow \sSet$---define the $n$-simplices $\Delta^n$ for $n \geq 0$, thought of as combinatorial avatars of the topological $n$-simplices. By the density of the Yoneda embedding, all other simplicial sets may be canonically presented as colimits of the simplices. There is a canonically-defined adjunction between simplicial sets and the strict 1-category of point-set topological spaces whose left adjoint ``geometrically realizes'' a simplicial set as a topological space:
\begin{equation}\label{eq:realization-sing}
\begin{tikzcd} & \DDelta \arrow[dr, "{\Delta^\bullet}"] \arrow[dl, hook', "\yo"'] \\ \sSet \arrow[rr, bend left=10, "{|-|\coloneq\colim_{\Delta^\bullet}}"] \arrow[rr, phantom, "\bot"] && \Top\,{.} \arrow[ll, bend left=10, "{\fun{Sing} \coloneqq \Hom(\Delta^\bullet, -)}"]
\end{tikzcd}
\end{equation}

There is something a bit strange about using simplicial sets as a combinatorial model for spaces, namely that the standard simplices are ``directed''; for instance there is no ``reversal'' map $\rho \colon \Delta^1 \to \Delta^1$ that exchanges the two vertices. While maps $p \colon \Delta^1 \to X$ are thought of as ``paths'' in $X$, paths cannot necessarily be reversed or composed unless $X$ is a \emph{Kan complex}, characterized by a right lifting property:
\[
    \begin{tikzcd} \Lambda^n_k \arrow[d, utcofarrow] \arrow[r, "x"] & X & & \Lambda^n_k \arrow[d, utcofarrow] \arrow[r, "x"] & X \arrow[d, two heads, "f"] \\ \Delta^n \arrow[ur, dashed] & & & \Delta^n \arrow[ur, dashed] \arrow[r, "y"'] & Y
    \end{tikzcd} \qquad \qquad n \geq 1, \quad 0 \leq k \leq n\,{.}
\]
More generally, a map of simplicial sets is called a \emph{Kan fibration} and denoted $f \colon X \fto Y$ if it has the right lifting property against all inclusions of \emph{horns} $\Lambda^n_k \subset \Delta^n$,  these being the simplicial subset spanned by all of the codimension 1 simplices except the one opposite vertex $k$. The Kan complexes define the so-called \emph{fibrant objects} and the Kan fibrations define the \emph{fibrations} in a Quillen model structure, that is, a strict 1-categorical presentation of the ($\infty$-)category of spaces. 

\begin{theorem}[{Quillen \cite{Quillen1967}}]\label{thm:quillen} There is a right proper simplicial model structure on simplicial sets whose cofibrations are the monomorphisms. 
\end{theorem}


\subsection{Families of spaces.}

The category of simplicial sets is \emph{locally cartesian closed}, meaning that for any $f \colon \Delta \to \Gamma$ the composition functor has two right adjoints defined by pullback along $f$ and pushforward along $f$:
\begin{equation}\label{eq:lcc}
\begin{tikzcd}[sep=large] \sSet_{/\Delta} \arrow[r, bend left=35, "\Sigma_f", "\bot"'] \arrow[r, bend right=35, "\Pi_f"', "\bot"] & \sSet_{/\Gamma}\,{.} \arrow[l, "f^*" description]
\end{tikzcd}
\end{equation}
Since $\sSet$ has a terminal object $\Delta^0$, this implies that $\sSet \cong \sSet_{/\Delta^0}$ and indeed any slice $\sSet_{/\Gamma}$ is cartesian closed and has all finite limits.

\begin{lemma} When $f \colon \Delta \fto \Gamma$ is a fibration, all three adjoints $\Sigma_f \dashv f^* \dashv \Pi_f$ preserve fibrations.\footnote{The fibration hypothesis is not needed for the pullback functor, but is needed for the other two adjoints.}
\end{lemma}

We introduce special notation for fibrations of simplicial sets $\rho \colon A \fto \Gamma$, which we think of as families of spaces parametrized by a base space $\Gamma$ referred to as the \emph{context}; note that the fibers will be fibrant objects but it not  necessary that the base $\Gamma$ is fibrant. Henceforth, we use ``family of spaces'' to refer to a fibration over a not necessarily fibrant base.

\begin{notation}\label{ntn:syntax} The following notation is intended to make the composition, pullback, and pushforward operations on families of spaces more legible.
  \begin{itemize}
    \item Fibrations, aka fibrant objects in a slice $\sSet_{/A}$, will be depicted as indexed families of spaces:
    \[ \begin{tikzcd} B \arrow[d, two heads, "{(B_a)_{a :A}}"]& & B_a \arrow[d, two heads] \arrow[dr, phantom, "\lrcorner" very near start] \arrow[r] & B \arrow[d, two heads, "{(B_a)_{a:A}}"] \\ A & & \Gamma \arrow[r, "a"'] & A\,{.}
    \end{tikzcd}
    \]
    Here, $B_a$ denotes the fiber over a generalized element $a \colon \Gamma \to A$, abbreviated $a:A$. We also refer to a family of spaces $(B_a)_{a :A}$ as a ``space $B$ in context $A$.''
    \item Pullbacks of families of spaces are also called \emph{substitutions}, as the fiber of the pullback of a family of spaces $(B_a)_{a :A}$ along $f \colon C \to A$  over a generalized element $c \colon \Gamma \to C$ is calculated by substituting the element $fc \colon \Gamma \to A$ for the element $a$:
    \[ \begin{tikzcd} B_f \arrow[r] \arrow[d, two heads, "{(B_{fc})_{c :C}}"'] \arrow[dr, phantom, "\lrcorner" very near start] & B \arrow[d, two heads, "{(B_a)_{a :A}}"] \\ C \arrow[r, "f"']& A \,{.}
    \end{tikzcd}
    \] 
     \item When $A$ is fibrant, the functor $\Sigma_A$ defined by composing with $! \colon A \fto \Delta^0$ sends a family of spaces $(B_a)_{a:A}$ to a space $\Sigma_{a : A} B_a$, namely the domain $B$ of the fibration:
     \[ \begin{tikzcd} \Sigma_{a:A}B_a \arrow[d, two heads, "{(B_a)_{a :A}}"'] \arrow[dr, phantom, "\mapsto"] & \Sigma_{a:A}B_a \arrow[d, two heads, "!"] \\ A  \arrow[r, two heads, "!"'] & \Delta^0\,{.}
    \end{tikzcd}
    \]
    More generally, when $A$ is a space $(A_\gamma)_{\gamma : \Gamma}$ in context $\Gamma$, the composition functor sends a fibration $(B_a)_{a :A}$ to the fibration:
    \[ \begin{tikzcd} \Sigma_{a:A}B_a \arrow[d, two heads, "{(B_a)_{a :A}}"'] \arrow[dr, phantom, "\mapsto"] & \Sigma_{\gamma : \Gamma}\Sigma_{a:A_\gamma}B_a \arrow[d, two heads, "{(\Sigma_{a: A_\gamma}B_a)_{\gamma : \Gamma}}"] \\ A  \arrow[r, two heads, "{(A_\gamma)_{\gamma :\Gamma}}"'] & \Gamma\,{.}
    \end{tikzcd}
    \]
    Since the effect of composing a pair of fibrations is to ``sum up the fibers'' we also refer to the composition functor as \emph{summation}.
  \item When $A$ is fibrant, the functor $\Pi_A$ defined by pushforward along $! \colon A \fto \Delta^0$ sends a fibration $(B_a)_{a:A}$ to a space $\Pi_{a : A} B_a$, namely the space $\Map_A(A,B)$ of \emph{sections} of the fibration:
  \[ \begin{tikzcd} \Sigma_{a:A}B_a \arrow[d, two heads, "{(B_a)_{a :A}}"'] \arrow[dr, phantom, "\mapsto"] & \Pi_{a:A}B_a \arrow[d, two heads, "!"] \\ A  \arrow[r, two heads, "!"'] & \Delta^0\,{.}
  \end{tikzcd}
  \]
  More generally, when $A$ is a space $(A_\gamma)_{\gamma : \Gamma}$ in context $\Gamma$, the pushforward functor sends a fibration $(B_a)_{a :A}$ to the fibration:
  \[ \begin{tikzcd} \Sigma_{a:A}B_a \arrow[d, two heads, "{(B_a)_{a :A}}"'] \arrow[dr, phantom, "\mapsto"] & \Sigma_{\gamma : \Gamma}\Pi_{a:A_\gamma}B_a \arrow[d, two heads, "{(\Pi_{a: A_\gamma}B_a)_{\gamma : \Gamma}}"] \\ A  \arrow[r, two heads, "{(A_\gamma)_{\gamma :\Gamma}}"'] & \Gamma\,{.}
  \end{tikzcd}
  \]
  \end{itemize}
\end{notation}

The notations just introduced are ``stable under substitution,'' commuting with pullbacks. By the \emph{Beck--Chevalley isomorphisms}, the fibers over $\gamma$ of the composition or pushforward along a fibration $(A_\gamma)_{\gamma : \Gamma}$ agrees with the spaces obtained by composition or pushforward over the space $A_\gamma$ defined by the pullback:
\[\begin{tikzcd}
    A_\gamma \arrow[d, two heads] \arrow[r] \arrow[dr, phantom, "\lrcorner" very near start] & A \arrow[d, two heads, "{(A_\gamma)_{\gamma : \Gamma}}"] \\ \Delta \arrow[r, "\gamma"'] & \Gamma\,{.}
  \end{tikzcd}
  \] 

\begin{lemma}\label{lem:beck-chevalley}
  The fibers  over a generalized element $\gamma \colon \Delta \to \Gamma$ of a sum or pushforward are calculated by:
\[ \begin{tikzcd} \Sigma_{a:A_\gamma}B_a \arrow[d, two heads] \arrow[r] \arrow[dr, phantom, "\lrcorner" very near start] & \Sigma_{\gamma : \Gamma}\Sigma_{a:A_\gamma}B_a \arrow[d, two heads, "{(\Sigma_{a: A_\gamma}B_a)_{\gamma : \Gamma}}"] & & \arrow[d, phantom, "\textup{and}"] & & \Pi_{a:A_\gamma}B_a \arrow[d, two heads] \arrow[r] \arrow[dr, phantom, "\lrcorner" very near start] & \Sigma_{\gamma : \Gamma}\Pi_{a:A_\gamma}B_a \arrow[d, two heads, "{(\Pi_{a: A_\gamma}B_a)_{\gamma : \Gamma}}"] \\ \Delta  \arrow[r, "{\gamma}"'] & \Gamma & & ~& & \Delta  \arrow[r, "{\gamma}"'] & \Gamma\,{.}
\end{tikzcd}
\]
\end{lemma}

\section{Path induction.}\label{sec:path}

We work with the spaces and families of spaces introduced above, with technical details made precise in reference to the structures defined in \cref{thm:quillen} and elaborated further in \cite{Riehl2024its,Shulman2015er,Shulman2019}.

For any space $A$, there is a natural \emph{path space factorization} of the diagonal map defined by exponentiation with the simplicial interval
\begin{equation}\label{eq:path-space}
    \begin{tikzcd}[sep=large] \Delta^0+\Delta^0 \arrow[r, tail] & \Delta^1 \arrow[r, we] & \Delta^0 & \rightsquigarrow & A \arrow[r, utcofarrow, "\refl"] & A^{\Delta^1} \arrow[r, two heads, "{(\ev_0,\ev_1)}"] & A \times A\,{.}
    \end{tikzcd}
  \end{equation}
This constructs a family of spaces over $A \times A$ whose total space is the path space $A^{\Delta^1}$ and whose fiber over $x, y : A$ we denote by $x \sim y$ as we think of this family of spaces as defining a generalized binary relation on the base space $A$. An element $p : x \sim y$ of the fiber is a \emph{path} from $x$ to $y$ in $A$
\[ 
  \begin{tikzcd} {x \sim y} \arrow[r]
    \arrow[d, two heads] \arrow[dr, phantom, "\lrcorner" very near start] & A^{\Delta^1} \arrow[d, two heads, "{(x \sim y)_{x, y : A}}"] \\ \Delta^0 \arrow[r, "{(x,y)}"'] & A \times A\,{.}
\end{tikzcd}
\]

As a consequence of \cref{thm:quillen}, the left map of \cref{eq:path-space}, the inclusion of constant paths, is a trivial cofibration in Quillen's model structure, and as such has the right lifting property with respect to any family of spaces:
\[
    \begin{tikzcd} A \arrow[d, tcofarrow, "\refl"'] \arrow[r, "e"] & E \arrow[d, two heads, "{\rho = (E_b)_{b : B}}"] & &\arrow[d, phantom, "\leftrightsquigarrow"] & \Map(A^{\Delta^1}, E) \arrow[d, utfibarrow, "\rho \circ - \circ \refl"] \\ A^{\Delta^1} \arrow[r, "b"'] \arrow[ur, dashed] \arrow[ur, dashed] & B &&  ~ & \Map(A,E) \times_{\Map(A,B)} \Map(A^{\Delta^1},B)\,{.}  \arrow[u, bend left, dashed]
      \end{tikzcd}        
\]
As a further consequence of \cref{thm:quillen}, this lifting property can be \emph{internalized}. The canonical map from the space $\Map(A^{\Delta^1},E)$ of maps from the path space of $A$ to $E$ to the space $\Map(A,E) \times_{\Map(A,B)} \Map(A^{\Delta^1},B) $ of commutative squares from $\refl$ to $\rho$ is a \emph{trivial fibration}, i.e., a ``contractible family of spaces.'' A section to this trivial fibration gives a continuous choice of solutions to lifting problems from $\refl$ to $\rho$. The fiber over a vertex $(e,b)$ defines a contractible space; thus solutions to lifting problems from $\refl$ to $\rho$ are ``homotopically unique.''

By pulling back along the codomain $b$, to solve lifting problems of the form above it suffices to consider lifting problems against families of spaces over the path space $A^{\Delta^1}$ whose codomain is the identity map:
\begin{equation}\label{eq:path-ind-as-section}
    \begin{tikzcd} A \arrow[rr, bend left, "e"]\arrow[d, tcofarrow, "\refl"'] \arrow[r, "d"] & P \arrow[d, two heads, "\chi"] \arrow[r] \arrow[dr, phantom, "\lrcorner" very near start] & E \arrow[d, two heads, "\rho"] & [-.75em] \arrow[d, phantom, "\leftrightsquigarrow"]   \\  A^{\Delta^1} \arrow[r, equals] \arrow[ur, dashed] & A^{\Delta^1} \arrow[r, "b"'] & B & ~   \end{tikzcd}   \quad
    \begin{tikzcd} \Map_{A^{\Delta^1}}(A^{\Delta^1},P) \arrow[d, utfibarrow, "- \circ \refl"] \arrow[r] \arrow[dr, phantom, "\lrcorner" very near start] & \Map(A^{\Delta^1},P) \arrow[d, utfibarrow, "\chi \circ - \circ \refl"] \\      \Map_A(A, P_\refl) \arrow[r]   \arrow[u, dashed, bend left, "\textup{path-ind}"] \arrow[d, two heads] \arrow[dr, phantom, "\lrcorner" very near start] & \Map(A,P) \times_{\Map(A,A^{\Delta^1})} \Map(A^{\Delta^1}, A^{\Delta^1})  \arrow[u, dashed, bend left] \arrow[d, two heads, "\pi"] \\ \Delta^0 \arrow[r, "\id"'] & \Map(A^{\Delta^1}, A^{\Delta^1})\,{.}
    \end{tikzcd}
  \end{equation}
Here the data of the lifting problem is given by a single map $d \colon A \to P$ that is a partial section of the fibration $\chi \colon P \fto A^{\Delta^1}$ over the map $\refl$. 
\[
\begin{tikzcd}
  & P \arrow[d, two heads, "{(P_{x,y,p})_{x,y:A, p: x \sim y}}"] && &\arrow[d, phantom, "\leftrightsquigarrow"] & [-.35em] P_\refl \arrow[r] \arrow[d, two heads, "{(P_{a,a,\refl_a})_{a:A}}"] \arrow[dr, phantom, "\lrcorner" very near start] & [+.65em] P \arrow[d, two heads, "{(P_{x,y,p})_{x,y:A, p: x \sim y}}"] && &\arrow[d, phantom, "\leftrightsquigarrow"] &[-.35em] A \arrow[r, "d", dashed] \arrow[d, tcofarrow, "\refl"'] & P \arrow[d, two heads, "{(P_{x,y,p})_{x,y:A, p: x \sim y}}"] \\ A \arrow[ur, "d", dashed] \arrow[r, "\refl"', tcofarrow] & A^{\Delta^1} & &&~ & A \arrow[u, bend left, "d", dashed] \arrow[r, "\refl"', tcofarrow] & A^{\Delta^1} & &&~ & A^{\Delta^1} \arrow[r, equals]  & A^{\Delta^1}\,{.}
\end{tikzcd}  
\]
As before, there exists a continuous choice of solutions to this lifting problem given by a section to the trivial fibration displayed above-center in \cref{eq:path-ind-as-section}. We call this map \emph{path induction}, using an analogy first observed by Awodey and Warren \cite{AwodeyWarren2009} and Gambino and Garner \cite{GambinoGarner2008}.  The existence of this map proves the following proposition:

\begin{proposition}[path induction, preliminary form]\label{prop:basic-path-ind} To define a section of a fibration over a path space $A^{\Delta^1}$, it suffices to define a partial section over the subspace $\refl \colon A \cwto A^{\Delta^1}$ of constant paths:
 \[
  \begin{tikzcd}[column sep=huge] A \arrow[r, "d"] \arrow[d, tcofarrow, "\refl"'] & P \arrow[d, two heads, "{(P_{x,y,p})_{x,y : A, p : x \sim y}}"] \\ A^{\Delta^1} \arrow[r, equals] \arrow[ur, dashed, "{\pathind(\refl \mapsto d)}" description] & A^{\Delta^1}\,{.}
  \end{tikzcd}
 \]
\end{proposition}

Path induction is a powerful proof technique, closely analogous to the principle of mathematical induction over the natural numbers. 
 By analogy with that case, it also has constructive content, allowing the ``recursive'' definition of morphisms, though we follow convention and use the term ``induction'' for both induction and recursion. 

\begin{construction}[inversion] Even though the simplicial interval $\Delta^1$ has no symmetries, when $A$ is a space, a path $p : x \sim y$ can be inverted to define a path $p^{-1} \colon y \sim x$. In other words, the path relation on spaces is symmetric as well as reflexive. This symmetry can be defined as a continuous function on path spaces by path induction, whence it suffices to specify the inverse of constant paths, which we take to be constant paths:
  \[ \begin{tikzcd} A \arrow[d, tcofarrow, "\refl"'] \arrow[r, utcofarrow, "\refl"] & A^{\Delta^1} \arrow[d, two heads, "{(y \sim x)_{x,y : A}}"]  \\ A^{\Delta^1} \arrow[r, two heads, "{(x \sim y)_{x,y : A}}"'] \arrow[ur, dashed, "{(-)^{-1}}" description] & A \times A
  \end{tikzcd}\qquad (-)^{-1} \coloneqq \pathind(\refl\mapsto\refl)\,{.}
  \]
\end{construction}

The full principle of path induction extends the version stated in \cref{prop:basic-path-ind} along two axes. The first extension uses the fact that Quillen's model structure satisfies the \emph{Frobenius condition} \cite{GambinoSattler2017}: The trivial cofibrations are stable under pullback along fibrations. The operation of pullback along a fibration can be thought of as introducing a trivial dependency and is referred to as \emph{weakening}. A more general version of path induction states that given a fibration over a weakening of a path space, to define a section it suffices to define a section over the subspace of constant paths. 

\begin{proposition}[path induction, intermediate form]\label{prop:weakened-path-ind} Given a composable pair of fibrations $\chi \colon Q \fto P$ and $\rho \colon P \fto A^{\Delta^1}$ over a path space, to define a section of $\chi$ it suffices to define a partial section over the pullback along $\rho$ of the subspace $\refl \colon A \cwto A^{\Delta^1}$ of constant paths:
  \[
  \begin{tikzcd}[column sep=huge]
 A \arrow[d, tcofarrow, "\refl"'] & \arrow[l, two heads, "\rho_\refl"'] \arrow[dl, phantom, "\llcorner" very near start] P_\refl \arrow[d, "{\refl_\rho}"' pos=.65, utcofarrow]  \arrow[r, "d"] & Q \arrow[d, two heads, "\chi = (Q_p)_{p : P}"] \\ 
    A^{\Delta^1} & \arrow[l, two heads, "\rho"] P \arrow[r, equals] \arrow[ur, dashed, "{\pathind(\refl_{\rho}\mapsto d)}" description] & P\,{.}
  \end{tikzcd}  
  \]
 \end{proposition}

\begin{construction}[concatenation] In a space $A$, paths $p : x \sim y$ and $q : y \sim z$ can be composed to define a path $p \ast q : x \sim z$. The composition function, establishing the transitivity of the path space relation, can be defined by path induction using the weakening of the reflexivity map displayed below-left:
  \[
  \begin{tikzcd}[column sep=6em]
 A \arrow[d, tcofarrow, "\refl"'] & \arrow[l, "{(\Sigma_{x:A}x \sim y)_{y : A}}"', two heads] A^{\Delta^1} \arrow[dl, phantom, "\llcorner" very near start] \arrow[d, "{(1, \refl)}"' pos=.65, tcofarrow]  \arrow[r, equals] & A^{\Delta^1} \arrow[d, two heads, "{(x \sim z)_{x,z : A}}"] \\  A^{\Delta^1} & \arrow[l, two heads, "{(\Sigma_{x:A}x \sim y)_{y,z : A, q : y \sim z}}"] A^{\Lambda^2_1} \arrow[r, two heads, "{(\Sigma_{y:A}x \sim y \sim z)_{x,z:A}}"'] \arrow[ur, dashed, "\ast"'] & A \times A
  \end{tikzcd}  \qquad \ast \coloneqq \pathind(p,\refl\mapsto p)\,{.}
  \]
  In summary, by path induction, the concatenation operation can be defined by specifying how to concatenate an arbitrary path $p : x \sim y$ by a constant path $\refl_y : y \sim y$, and we define concatenation with a constant path to be the identity function.
\end{construction}

This construction of the concatenation function avoids the use of any higher simplices, though by homotopical uniqueness of solutions to lifting problems it is equivalent to the map defined by the more standard construction:
\begin{equation}\label{eq:concat-alt}
  \begin{tikzcd} \Lambda^2_1 \arrow[r, tcofarrow] & \Delta^2 & \Delta^1 \arrow[l, tail, "\delta^1"'] & & \rightsquigarrow & & \ast \colon  A^{\Lambda^2_1} \arrow[r, tcofdashed, bend left] & A^{\Delta^2} \arrow[l, tfibarrow] \arrow[r, two heads, "-\cdot\delta^1"] & A^{\Delta^1} \,{.} \end{tikzcd}
\end{equation}

\begin{construction}[transport]\label{cons:tr}
Families of spaces $(B_a)_{a : A}$ are accompanied by path lifting or ``transport'' operations that lift a path $p : x \sim y$ in the base space $A$ to a continuous map $\tr_p \colon B_x \to B_y$ between the fibers. 

As before, this operation can be expressed as a continuous function between spaces, where the domain is formed by pulling the domain projection $\ev_0 \colon A^{\Delta^1} \fto A$ back along $\rho \colon B \fto A$:
\[
  \begin{tikzcd}[column sep=6em]
 A \arrow[d, tcofarrow, "\refl"'] & \arrow[l, "{(B_a)_{a :A}}"', two heads] B \arrow[dl, phantom, "\llcorner" very near start] \arrow[d, "{\refl \times 1}"' pos=.65, tcofarrow]  \arrow[r, equals] & B \arrow[d, two heads, "{\rho = (B_y)_{y :A}}"] \\   A^{\Delta^1}  & \arrow[l, two heads, "{(B_x)_{x,y : A, p: x \sim y}}"] A^{\Delta^1}\!\!\times_A\! B \arrow[r, two heads, "{(\Sigma_{x:A} (x \sim y) \times B_x)_{y : A}}"'] \arrow[ur, dashed, "\tr"'] & A
  \end{tikzcd}  \qquad \tr \coloneqq \pathind(\refl,u \mapsto u)\,{.}
  \]
In summary, the transport operation is defined by path induction by declaring that transport along a constant path is the identity function.
\end{construction}

The second extension of the homotopy type theoretic principle of path induction allows $A$ to be a space in context $\Gamma$, starting from an arbitrary family of spaces $(A_\gamma)_{\gamma : \Gamma}$, whether or not the base $\Gamma$ is fibrant. In context $\Gamma$, we replace the path space factorization \cref{eq:path-space} by the fibered path space factorization, a construction which is stable under substitution:
\begin{equation}\label{eq:relative-path-space-factorization}
    \begin{tikzcd}[sep=tiny]
        A \arrow[rrr] \arrow[dddd, two heads] \arrow[dr, tcofdashed] & & & A^{\Delta^1} \arrow[rrr] \arrow[dddd, two heads] & & & A \times A \arrow[dddd, two heads] \\ &  \Delta^1 \pitchfork_\Gamma A \arrow[dddl, two heads, dotted] \arrow[dddr, phantom, "\lrcorner" very near start]\arrow[urr, dotted] \arrow[dr, two heads, dashed ]
        \\ & & A \times_\Gamma A \arrow[ddll, two heads, dotted] \arrow[uurrrr, dotted]  \arrow[ddr, phantom, "\lrcorner" very near start]\\ ~
        \\ \Gamma \arrow[rrr] & & ~ & \Gamma^{\Delta^1} \arrow[rrr] & & & \Gamma \times \Gamma
    \end{tikzcd} \qquad 
        \begin{tikzcd}[sep=tiny]  A_f \arrow[dr, tcofarrow] \arrow[ddddrr, phantom, "\lrcorner" pos=.05] \arrow[ddddr, two heads] \arrow[rrr, dotted] & & &[-35pt] A \arrow[dr, tcofarrow] \arrow[ddddr, two heads]\\ & \Delta^1 \pitchfork_\Delta A_f \arrow[dr, two heads] \arrow[ddd, two heads] \arrow[rrr,dotted] \arrow[dddr, phantom, "\lrcorner" very near start] & & & \Delta^1 \pitchfork_\Gamma A \arrow[dr, two heads] \arrow[ddd, two heads]\\ & & A_f \times_\Delta A_f \arrow[ddl, two heads] \arrow[rrr, dotted] \arrow[ddr, phantom, "\lrcorner" very near start]& & & A \times_\Gamma A \arrow[ddl, two heads] \\ ~ \\ & \Delta \arrow[rrr, dotted, "f"'] &~ & ~& \Gamma\,{.}
    \end{tikzcd} 
  \end{equation}

The inclusion $\refl \colon A \cwto \Delta^1 \pitchfork_\Gamma A$ is a trivial cofibration so we have a lifting property exactly as above. This proves the general form of the path induction principle.

\begin{proposition}[path induction, final form]\label{prop:path-ind}
  Given a composable pair of fibrations $\chi \colon Q \fto P$ and $\rho \colon P \fto \Delta^1 \pitchfork_\Gamma A$ over a path space in any context $\Gamma$, to define a section of $\chi$ it suffices to define a partial section over the pullback along $\rho$ of the subspace $\refl \colon A \cwto \Delta^1 \pitchfork_\Gamma A$ of constant paths:
  \[
  \begin{tikzcd}[column sep=huge]
    P_\refl \arrow[d, "{\refl_\rho}"', tcofarrow] \arrow[r, "\rho_\refl", two heads] \arrow[dr, phantom, "\lrcorner" very near start] & A \arrow[d, utcofarrow, "\refl"] & \arrow[d, phantom, "\rightsquigarrow"] & P_\refl \arrow[d, "{\refl_\rho}"', tcofarrow]  \arrow[r, "d"] & [+1em] Q \arrow[d, two heads, "\chi = (Q_p)_{p : P}"] \\ 
    P \arrow[r, two heads, "\rho"'] & {\Delta^1}\pitchfork_\Gamma A & ~ & P \arrow[r, equals] \arrow[ur, dashed, "{\pathind(\refl_{\rho}\mapsto d)}" description] & P\,{.}
  \end{tikzcd}  
  \]
\end{proposition}

\begin{remark}\label{rmk:stable-under-slicing}
 \cref{prop:path-ind} is a precise analog of \cref{prop:weakened-path-ind} just interpreted in the slice 1-category $\sSet_{/\Gamma}$ instead of in $\sSet$.  Crucially, all of the properties of \cref{thm:quillen} are stable under slicing, meaning inherited by all sliced 1-categories $\sSet_{/\Gamma}$. In particular, this is why we described Quillen's model structure as \emph{simplicial} (enriched over Quillen's model structure on simplicial sets) rather than \emph{cartesian closed} (enriched over itself), because the former property is stable under slicing while the latter is not.
\end{remark}

\section{Univalent universes.}\label{sec:ua}

In \S\ref{ssec:props}, we construct spaces whose points encode proofs of mathematical propositions, following \cite[\S4]{Shulman2015er}, where more details can be found. In \S\ref{ssec:uu}, we introduce Voevodsky's univalence axiom and sketch a modern approach to constructing univalent universes due to Shulman \cite{Shulman2019}. In \S\ref{ssec:ua-applications}, we tour a few standard applications, which are described in more detail in the textbooks \cite{hottbook, Rijke2025,Symmetry}.

\subsection{Propositions as spaces.}\label{ssec:props}

The adjoint triples $\Sigma_f \dashv f^* \dashv \Pi_f$ \eqref{eq:lcc} can be used to construct spaces in the empty context or fibrations in a general context $\Gamma$ whose global elements---points in the former case and sections in the latter case---provide data witnessing the proof of some proposition. Thus, these spaces exhibit mathematical propositions as homotopy types. 
By \cref{rmk:stable-under-slicing}, our constructions equally make sense in slices $\sSet_{/\Gamma}$, so  we leave the context $\Gamma$ implicit in our notation, representing a family of spaces $(A_\gamma)_{\gamma : \Gamma}$ as a space ``$A$'' and reverting to the notation $A^{\Delta^1}$ for the path space in context $\Gamma$. 

\begin{construction}[contractibility]\label{cons:iscontr} Fix a space $A$.
  The pushforward of the path space family $(x\sim y)_{x, y :A}$ along the projection away from the second coordinate defines the family of spaces $(\Pi_{y : A}x\sim y)_{x : A}$, and then the sum along the projection away from the first coordinate defines the space $\Sigma_{x :A} \Pi_{y :A} x\sim y$.
\[ \begin{tikzcd} A^{\Delta^1} \arrow[d, two heads, "{(x\sim y)_{x,y :A}}"] \arrow[dr, phantom, "\mapsto"] & [+2em]  \Sigma_{x :A} \Pi_{y : A} x \sim y \arrow[d, two heads, "{(\Pi_{y:A}x \sim y)_{x : A}}"] \arrow[dr, phantom, "\mapsto"] & [+2em] \Sigma_{x: A} \Pi_{y :A} x \sim y \arrow[d, two heads, "!"]\\ A \times A & A & \Delta^0\,{.}
\end{tikzcd}
\]
A global element of $\Sigma_{x:A} \Pi_{y :A} x \sim y$ defines an element $a : A$ together with an element of the fiber $\Pi_{y :A} a \sim y$, which is a section of the based path space fibration. This data proves the contractibility of $A$, and thus we define
  \[  \isContr(A) \coloneqq \Sigma_{x:A} \Pi_{y :A} x \sim y \, .\]
\end{construction}

An example of a contractible space is a based path space, defined for a fixed point $a : A$ by 
\[  \begin{tikzcd}
    \Sigma_{y:A} a \sim y \arrow[d, two heads] \arrow[r] \arrow[dr, phantom, "\lrcorner" very near start] & A^{\Delta^1} \arrow[d, two heads, "{(\Sigma_{y : A} x \sim y)_{x : A}}"] \\ \Delta^0 \arrow[r, "a"'] & A\,{.}
  \end{tikzcd}
\]
A proof that $\Sigma_{y:A} a \sim y$ is contractible defines a global element of $\isContr(\Sigma_{y:A} a \sim y)$, but we can do better. When we ``fix a point $a : A$'' or ``let $a$ be a point $A$,'' we are introducing the space $A$ as the context. Thus, our based path spaces are more properly thought of as families of spaces  $(\Sigma_{y:A} a \sim y)_{a : A} \in \sSet_{/A}$.

By \cref{rmk:stable-under-slicing}, the construction of \cref{cons:iscontr} can be interpreted in any context. When applied to a fibration $(A_\gamma)_{\gamma : \Gamma}$ it defines a fibration $\isContr((A_\gamma)_{\gamma : \Gamma}) \cong (\isContr(A_\gamma))_{\gamma : \Gamma} \in \sSet_{/\Gamma}$ whose fiber over $\gamma : \Gamma$ is the space $\isContr(A_\gamma)$ by \cref{lem:beck-chevalley}. When applied to the fibration $(\Sigma_{y:A} a \sim y)_{a : A}$, we obtain a fibration $(\isContr(\Sigma_{x :A} a \sim x))_{a : A} \in \sSet_{/A}$ admitting a section, which proves the contractibility of the based path spaces, simultaneously and continuously for all $a :A$. Indeed:

\begin{lemma}\label{lem:iscontr} A family of spaces $(A_\gamma)_{\gamma : \Gamma}$ defines a contractible family of spaces if and only if $(\isContr(A_\gamma))_{\gamma : \Gamma} \in \sSet_{/\Gamma}$ has a section.
\end{lemma}

\begin{construction}[equivalence]
  Consider a map $f \colon A \to B$ between spaces. The familiar mapping path space construction below-left defines a family of spaces $(\fib_fb)_{b : B}$ as below-right:
\[     \begin{tikzcd}
    \Sigma_{a :A}\Sigma_{b :B}fa \sim b \arrow[r] \arrow[d, two heads, "{(fa \sim b)_{a : A, b:B}}"'] \arrow[dr, phantom, "\lrcorner" very near start] & B^{\Delta^1} \arrow[d, two heads, "{(x\sim y)_{x,y : B}}"] & & \arrow[d, phantom, "\rightsquigarrow"] & \fib_fb \coloneqq \Sigma_{a :A} fa \sim b \arrow[r] \arrow[d, two heads] \arrow[dr, phantom, "\lrcorner" very near start ] &  \Sigma_{b :B} \Sigma_{a :A}fa \sim b \arrow[d, two heads, "{(\Sigma_{a:A} fa \sim b)_{b : B}}"]  \\ A \times B \arrow[r, "f \times \id"'] & B \times B & & ~ & \Delta^0 \arrow[r, "b"'] & B
  \end{tikzcd}
\]
This allows us to define the space
  \[ \isEquiv(f) \coloneqq \Pi_{b : B} \isContr(\fib_fb).\]
\end{construction}

As the terminology suggests, a point in the space $\isEquiv(f)$ proves that $f \colon A \to B$ is an equivalence between spaces. Such a point provides the data of a section to the fibration below-left:
\[
  \begin{tikzcd}
    \Sigma_{b :B} \isContr(\fib_fb) \arrow[d, two heads, "{(\isContr(\fib_fb))_{b:B}}"] & & \arrow[d, phantom, "\rightsquigarrow"] &    \Sigma_{b :B} \fib_fb \arrow[d, two heads, "{(\fib_fb)_{b:B}}"] \\ B \arrow[u, bend left, dashed] & & ~ & \arrow[u, dashed, bend left] B\,{.}
  \end{tikzcd}
  \] 
  Passing to the center of contraction, this defines a section to the fibration above-right, which gives the data of a continuous function $g \colon B \to A$ together with a homotopy $\beta : \Pi_{b : B} fg b \sim b$. The remaining data---which witnesses that the fibers are contractible, not just inhabited---can be used to construct a second homotopy $\alpha: \Pi_{a : A} gfa \sim a$. 

  Again, this construction can be interpreted in any context:

  \begin{lemma}\label{lem:isequiv}
    A map $f$ between families of spaces over $\Gamma$ is an equivalence if and only if $\isEquiv(f) \in \sSet_{/\Gamma}$ has a section.
  \end{lemma}

  \begin{construction}[equivalences]\label{cons:eq} For spaces $A$ and $B$, define a space
  \[ A \simeq B \coloneqq \Sigma_{f : A \to B} \isEquiv(f).\]
Here, we are using $A \to B$ as alternate notation for the mapping space $\Map(A,B) \cong B^A$. The space $A \simeq B$ has a natural projection map to the mapping space $A \to B$, which is a homotopy monomorphism because the space $\isEquiv(f)$ is a \emph{proposition}: if $\isEquiv(f)$ contains any points, then it is contractible. 
\end{construction}

\subsection{Constructing univalent universes.}\label{ssec:uu}

The previous constructions were defined in the context of arbitrary spaces $A$ and $B$. These can be regarded as constructions in context $\univ \times \univ$,  where we use a \emph{universe} $\univ$ of small\footnote{Here ``small'' means that the cardinality of the fibers is bounded by some infinite regular cardinal $\kappa$---an implicit parameter in everything that follows---which is often assumed to be inaccessible so that the universe is closed under all of the desired type theoretic operations. See \cite[\S 4]{Shulman2019} for an abstract treatment of the cardinality restrictions required to define universes.
}  spaces to define our context. Our next task is to define this space $\univ$.

What we seek is both a space $\univ$ of small spaces and a universal family $(A)_{A : \univ}$ of small spaces so that generalized elements $A \colon \Gamma \to \univ$ define families of small spaces in context $\Gamma$: 
\begin{equation}\label{eq:univ-pullback}
  \begin{tikzcd}
    A \arrow[d, two heads, "{(A_\gamma)_{\gamma : \Gamma}}"'] \arrow[dr, phantom, "\lrcorner" very near start] \arrow[r] & \tilde{\univ} \arrow[d, two heads, "{(A)_{A : \univ}}"] \\ \Gamma \arrow[r, "A"'] & \univ\,{.}
  \end{tikzcd}
\end{equation}
The universal family of small spaces is known as an \emph{object classifier}, in analogy with the \emph{subobject classifier} that can be found in any elementary 1-topos. 

Following Voevodsky \cite{KapulkinLumsdaine2021}, the universe is the fundamental ingredient used to solve the coherence problems inherent in interpreting the strict syntax of homotopy type theory (where in particular, the substitution operation is defined up to \emph{equality}) in a 1-category, where constructions like the pullback \cref{eq:univ-pullback} are only defined up to isomorphism, or even worse in an ($\infty$-)category, where constructions like \eqref{eq:univ-pullback} are only well-defined up to equivalence. Hence, we construct the universe in a strict 1-categorical presentation of the ($\infty$-)category of spaces, such as given by \cref{thm:quillen}, and require the classifying pullbacks to be 1-categorical rather than ($\infty$-)categorical. 

Na\"{i}vely, it appears that the 1-categorical universal property of the universal family $(A)_{A : \univ}$ states that it represents small families of spaces. But this is not quite correct, because small families of spaces with the substitution action define a 1-groupoid valued pseudofunctor, rather than a strict set-valued functor; see \cite[\S 4.1]{Riehl2024its} for discussion. Shulman introduces an axiomatic \emph{notion of fibred structure} \cite[\S3]{Shulman2019}, which in particular defines a 1-groupoid-valued pseudofunctor $\mathbb{F}$. Examples include pullback stable families of morphisms, like the fibrations in a Quillen model category, but also examples where the notion of ``fibration'' is encoded by additional structure borne by, rather than a mere property of, a map. In a Grothendieck 1-topos with a Cisinski model structure like in \cref{thm:quillen}, any notion of fibred structure $\mathbb{F}$ that satisfies a size condition---is \emph{small groupoid valued}---and a locality condition---is \emph{locally representable}---admits a universe $\univ$ with a map $\yo_{\univ} \to \mathbb{F}$ that is a ``surjective equivalence'' in a suitable sense \cite[5.10]{Shulman2019}. By the Yoneda lemma, this data encodes a fibration (in the sense of the notion of fibred structure $\mathbb{F}$) over $\univ$, defining the desired universal family $\upsilon \colon \tilde{\univ} \fto \univ$ of spaces.

The univalence axiom characterizes the paths $(A \sim B)_{A, B : \univ}$ in the universe in terms of another universal family. Consider the family of spaces $(A \simeq B)_{A, B : \univ}$ encoding the universal equivalence, as defined by \cref{cons:eq}. By construction, a generalized element 
\[
  \begin{tikzcd}  A \simeq B \arrow[r] \arrow[dr, phantom, "\lrcorner" very near start] \arrow[d, two heads] & \Eq(\tilde{\univ}) \arrow[d, two heads, "{(A\simeq B)_{A,B : \univ}}"] \\ \Gamma \arrow[u, bend left, dashed, "{e}"] \arrow[r, "{(A,B)}"']& \univ \times \univ
  \end{tikzcd} \]
encodes an equivalence $e \colon A \simeq B$ between the families of spaces over $\Gamma$ classified by the maps $A, B \colon \Gamma \to \univ$.

\begin{definition}[{Voevodsky \cite{KapulkinLumsdaine2021}}]\label{defn:univalence}
  The fibration $\upsilon \colon \tilde{\univ} \to \univ$ is \textbf{univalent} if the comparison map from the path space to the space of equivalences defined by path induction is an equivalence:
\[
  \begin{tikzcd}[column sep=large] \univ\arrow[d, "\refl"', tcofarrow] \arrow[r, "\id"] & \Eq(\tilde{\univ}) \arrow[d, two heads,"{(A\simeq B)_{A,B: \univ}}"] \\ \univ^{\Delta^1} \arrow[ur, dashed, "\patheq" description] \arrow[r, two heads, "{(A \sim B)_{A, B: \univ}}"'] & \univ \times \univ
  \end{tikzcd} \qquad \patheq \coloneqq \pathind(\refl \mapsto \id)\,{.}
  \]
\end{definition}

\subsection{Applications of Univalence.}\label{ssec:ua-applications}

Univalence implies that the synthetic language for spaces is natively ``derived,'' with all constructions respecting equivalences between spaces. A construction on spaces can be interpreted as defining a family of spaces $(E_X)_{X: \univ}$ over the universe $\univ$. By \cref{cons:tr}, a path $p \colon X \sim Y$ between spaces $X,Y : \univ$ defines a transport function $\tr_p \colon E_X \to E_Y$. By another application of path induction, this transport function 
is an equivalence, since this is true in the case of transport along constant paths, given by the identity function. By univalence, an equivalence $X \simeq Y$ can be converted into a path $X \sim Y$ giving rise to an equivalence $E_X\simeq E_Y$.

The universe can also be used to define other classifying spaces.  Univalence is then used to characterize their path spaces  by means of a family of theorems known as the \emph{structure identity principle} \cite{CoquandDanielsson2013, AhrensNorthShulmanTsementzis2020}.

\begin{example}\label{ex:B-sigma-n}
  For any natural number $n$, 
define $\Fin_n \coloneqq \Sigma_{X : \univ} \| \mathbf{n} \simeq X \|$ to be the space of $n$ element spaces, i.e., the subuniverse of spaces that are equivalent to the discrete space $\mathbf{n}$ on $n$ points. Here, the ``$\| \|$'' denotes the \emph{propositional truncation}, meaning that the data of an explicit equivalence $\mathbf{n} \simeq X$ is not recorded. It follows that the natural projection map from $\Fin_n$ to $\univ$ is a homotopy monomorphism. The space $\Fin_n$ is naturally regarded as a based space, with $\mathbf{n} : \Fin_n$ serving as the base point.
  
  By the univalence axiom, $\Fin_n$ is a 1-type, with all homotopy groups vanishing above dimension 1. Paths  in $\Fin_n$ are equivalent to paths  in $\univ$ between the objects in $\Fin_n$; in particular, since all $X, Y : \Fin_n$ are equivalent to $\mathbf{n}$ and thus to each other, the space $\Fin_n$ is $0$-connected, i.e., has a single path component. By univalence, paths $X \sim Y$ in $\Fin_n$ are equivalent to equivalences $X \simeq Y$. It follows that the based loop space $\Omega(\Fin_n) \coloneqq \mathbf{n} \sim \mathbf{n}$ is equivalent to the space of equivalences $\mathbf{n} \simeq \mathbf{n}$, better known as the symmetric group $\Sigma_n$ on $n$ elements. 
\end{example}

As noted in \cref{ex:B-sigma-n}, the space $\Fin_n$ is a pointed, 0-connected, 1-type whose loop space is the discrete space $\Sigma_n$. Thus, we have the alternate notation $B\Sigma_n \coloneqq \Fin_n$ that is more amenable to generalization. 

\begin{example}\label{ex:BG}
  For a discrete group $G$, there is a standard ``analytic'' construction of its classifying space as the geometric realization of the nerve of the strict 1-category with one object whose endoarrows are elements of the group $G$. As this construction involves a left adjoint \cref{eq:realization-sing}, it is somewhat delicate to establish the mapping in universal property of the space so-constructed. 
Using the language for synthetic spaces, we can instead construct $BG$ as the space spanned by those $X : \univ$ that are $G$-torsors, those $0$-types $X$ equipped with a free and transitive $G$-action. To define a family of $G$-torsors over a base space $\Gamma$ is then just to define a map $X \colon \Gamma \to BG$. 
\end{example}

\begin{example}
  Extending the perspective on groups introduced in \cref{ex:BG}, we may define the space of groups to be the space of pointed, 0-connected, 1-types. Alternatively, we can define the space of groups following the standard definition: ``a group is a discrete space, with a binary multiplication, satisfying axioms'' as we illustrate in the simpler case of $H$-spaces:\footnote{The extra coherence $\lambda_e \sim \rho_e$ is added to the usual definition to give a better-behaved moduli space of $H$-spaces  \cite{BuchholtzChristensenFlatenRijke2025}.}
 \[ H\text{-}\type{Spc} \coloneqq \Sigma_{X : \univ} \Sigma_{e : X} \Sigma_{\mu : X \times X \to X} \Pi_{x : X} \Sigma_{\lambda_x : \mu(e,x) \sim x} \Sigma_{\rho_x :\mu(x,e) \sim x} \lambda_e \sim \rho_e .\] 
  By univalence, paths in $H\text{-}\type{Spc}$ are equivalences of $H$-spaces via $H$-space homomorphisms; this is an instance of the aforementioned structure identity principle. Similarly, paths in $\type{Group}$ are group isomorphisms \cite{Rijke2025}.

Both spaces of groups just described---the ``concrete'' groups, as encoded by their classifying spaces $BG$, and the ``abstract'' groups, as captured by the usual set-theoretic definition---are equivalent.   For a development of group theory from the perspective of classifying spaces, see \cite{Symmetry}. 
\end{example}

An advantage of the classifying space approach is that it allows many constructions of group theory to extend to \emph{higher groups}. If $a :A$ is a point in an arbitrary space $A$, we may form the space $B\Aut(a) \coloneqq \Sigma_{x : A} \| a \sim x\|$ by generalizing the construction used in \cref{ex:B-sigma-n}. This is again a pointed space, with basepoint $a : B\Aut(a)$, and one recovers the higher automorphism group $\Aut(a)$ as $\Omega(B\Aut(a)) \coloneqq a \sim a$. If $A$ is an $(n+1)$-type, then so is $B\Aut(a)$, and the based loop space $\Aut(a)$ is $n$-truncated. This constructs an example of an $(n+1)$-\emph{group}, where the notion of $1$-group recovers the standard examples of automorphism groups of discrete spaces. In particular, we may form $B\Aut(A)$ associated to $A : \univ$ and immediately from its definition a map $X \to B\Aut(A)$ is a family of spaces over $X$ whose fibers are equivalent to the space $A$.

  In general, we define a \emph{higher group} to be a pointed connected space $BG$, where $G\coloneqq \Omega BG$ is the carrier space for the higher group and $BG$ is the classifying space. Groups with untruncated classifying spaces are $\infty$-\emph{groups}, while groups with $(n+1)$-truncated classifying spaces are $(n+1)$-groups.  An \emph{action} of a higher group $BG$ on a space $X$ is just a map $X : BG \to \univ$ sending the basepoint of $BG$ to $X : \univ$. Thus an action provides a family of spaces $(X_z)_{z : BG}$ over $BG$. The \emph{invariants} and \emph{coinvariants} of this action are given by the spaces $X^{hG} \coloneqq \Pi_{z : BG} X_z$ and $X_{hG} \coloneqq X{\sslash}G \coloneqq \Sigma_{z :BG} X_z$ respectively. See \cite{BuchholtzDoornRijke2018} for considerably more. 

\section{A language for synthetic categories}\label{sec:cats}

In the language for synthetic categories---introduced in \cite{RiehlShulman2017} and now commonly referred to as \emph{simplicial homotopy type theory} or \emph{simplicial type theory}---we extend the language for synthetic spaces  with an axiomatic directed interval $\2$ satisfying strict interval axioms. That is $\2$ is equipped with distinct elements $0,1 : \2$ and a binary predicate $(x \leq y)_{x,y : \2}$ that is reflexive, antisymmetric, transitive, and total, with minimum 0, maximum 1. We think of the object $\2$ as the ``walking arrow,'' with source 0 and target 1. Following the construction of Joyal's proof that the classifying 1-topos of the theory of the strict interval is the 1-category of simplicial sets \cite[\S VIII.8]{MacLaneMoerdijk1994}, we build a ``walking composable pair of arrows'' $\3$---which we identify as the subspace of the square $\2 \times \2$ carved out by the predicate $(y \leq x)_{x,y : \2}$---or more generally, the ``walking composable sequence of $n$ arrows'' $\mathbbe{n}+\1$ for each $n \geq 0$. 
The language may also be used to construct \emph{face inclusions} $\delta^i \colon \mathbbe{n} \to \mathbbe{n}+\1$ and \emph{degeneracy maps} $\sigma^j \colon \mathbbe{n}+\1 \to \mathbbe{n}$ defining the cosimplicial object below-left:
\begin{equation}\label{eq:categorical-simplices} \begin{tikzcd} \1 \arrow[r, shift left=1em, "\delta^1" description] \arrow[r, shift right=1em, "\delta^0" description] & \2 \arrow[l, "\sigma^0" description] \arrow[r, "\delta^1" description] \arrow[r, shift right=2em, "\delta^0" description] \arrow[r, shift left=2em, "\delta^2" description] & \3 \arrow[l, shift right=1em, "\sigma^0" description] \arrow[l, shift left=1em, "\sigma^1" description] & \cdots & \eqqcolon & \Delta^0 \arrow[r, shift left=1em, "\delta^1" description] \arrow[r, shift right=1em, "\delta^0" description] & \Delta^1 \arrow[l, "\sigma^0" description] \arrow[r, "\delta^1" description] \arrow[r, shift right=2em, "\delta^0" description] \arrow[r, shift left=2em, "\delta^2" description] & \Delta^2 \arrow[l, shift right=1em, "\sigma^0" description] \arrow[l, shift left=1em, "\sigma^1" description] & \cdots \end{tikzcd}
\end{equation}

We may think of this simplicial structure as freely added to the very convenient category of spaces introduced in \S\ref{sec:spaces}, and thus our basic semantic setting is $\sSet^{\DDelta^\op}$---a strict 1-category whose objects we refer to as \emph{simplicial spaces}. This is again an $\infty$-topos, and in particular a locally cartesian closed category, for which all of the constructions \S\ref{sec:spaces}--\ref{sec:ua} apply.  In $\sSet^{\DDelta^\op}$
the walking arrow is interpreted by the discretely embedded representable $\Delta^1$ and similarly for the other objects in the diagram \cref{eq:categorical-simplices}. Indeed, it is conventional to reuse the simplicial notation for the walking arrows just constructed---$\Delta^0 \coloneqq \1$, $\Delta^1 \coloneqq \2$, $\Delta^2 \coloneqq \3$, etc---as indicated on the right-hand side of \cref{eq:categorical-simplices}.  More generally, we'll use simplicial notation such as $\Lambda^2_1\subset \Delta^2$ for discretely embedded simplicial spaces, i.e., for simplicial objects valued in discrete spaces. 

\begin{warning} Beware that the object used to define the path space factorization in \cref{eq:path-space} is distinct from the simplicial space $\Delta^1 \coloneqq \2$ used here: that object is a constant simplicial object valued in a non-discrete space. Because the map $\refl$ defining the inclusion of the subspace of constant paths is an equivalence, we may use $A$ as an equivalent representation of the total space of the path space family on $A$ on the occasions it appears here. 
\end{warning}

In an attempt to be as concrete as possible, we'll refer to ``simplicial spaces'' in what follows, but the semantics of the language of synthetic category theory is more general, interpretable in the category of simplicial objects in any $\infty$-topos \cite{RiehlShulman2017,Weinberger2022ss}.

\begin{construction}[arrows]\label{cons:homs}
  For a simplicial space $C$, we define a family of simplicial spaces $(\Hom(x,y))_{x,y : C}$ by exponentiation with the walking arrow: 
  \[ 
    \begin{tikzcd}[sep=large] \Delta^0+\Delta^0 \arrow[r, tail] & \Delta^1 \arrow[r] & \Delta^0 & \rightsquigarrow & C \arrow[r, "\id"] & C^{\Delta^1} \arrow[r, two heads, "{(\ev_0,\ev_1)}"] & C \times C\, {.}
    \end{tikzcd}
\] 
  We refer to an element $f : \Hom(x,y)$ of the fiber as an \emph{arrow} from $x$ to $y$ in $C$.
\[ 
  \begin{tikzcd} {\Hom(x,y)} \arrow[r]
    \arrow[d, two heads] \arrow[dr, phantom, "\lrcorner" very near start] & C^{\Delta^1} \arrow[d, two heads, "{(\Hom(x,y) )_{x, y : C}}"] \\ \Delta^0 \arrow[r, "{(x,y)}"'] & C \times C\, {.}
\end{tikzcd}
\]
\end{construction}

Thus, all simplicial spaces $C$ have families of arrows $(\Hom(x,y))_{x,y : C}$ as well as identity arrows $\id_x : \Hom(x,x)$ for each $x : C$ defined by the partial section $\id \colon C \to C^{\Delta^1}$. But not all simplicial spaces should be thought of as categories. A \emph{precategory} is a simplicial space $C$ so that all pairs of composable arrows have a unique composite as is made precise by the following definition.

\begin{definition}\label{defn:segal}
A simplicial space $C$ is a \textbf{precategory} if the map $C^{\Delta^2} \fto C^{\Lambda^2_1}$ defined by the inclusion $\Lambda^2_1 \hookrightarrow \Delta^2$ is an equivalence.  
\end{definition}

  A precategory $C$ has a canonical composition function defined for any triple of elements $x,y,z : C$ 
  \[ \begin{tikzcd} \Hom(y,z) \times \Hom(x,y) \arrow[r, "\circ"] & \Hom(x,z)\, {.}\end{tikzcd}\] 
  This map is defined by identifying a composable pair of arrows $g,f$ with an element $(g,f) : C^{\Lambda^2_1}$, then using the inverse of the equivalence of \cref{defn:segal} to extend this data to a pair of composable arrows $\comp(g,f) : C^{\Delta^2}$, and finally precomposing with the function $\delta^1 \colon \Delta^1 \to \Delta^2$ to obtain an arrow $g \circ f : C^{\Delta^1}$ from $x$ to $z$, i.e., by the analogous construction of \cref{eq:concat-alt}.  But there are several advantages to giving the definition in the form that we did.  
  
  \begin{lemma}
  In a precategory $C$:
  \begin{itemize}
    \item For all $f : \Hom(x,y)$, $f \circ \id_x \sim f$ and  $\id_y \circ f \sim f$.
    \item For all composable triples of arrows $h,g,f$, $(h \circ g) \circ f \sim h \circ (g \circ f)$.
  \end{itemize}  
  \end{lemma} 
  
  These properties both follow from the higher structure inherent in the equivalence of simplicial spaces of \cref{defn:segal} and would not follow automatically from an axiomatized binary composition function.  Another useful consequence is ``functoriality for free,'' essentially by commutativity of pre- and post-composition in the category of simplicial spaces:

  \begin{lemma}\label{lem:fun-for-free}
    Any function $F \colon C \to D$ between precategories is automatically functorial, preserving both identities and composition.
  \end{lemma}

The simplicial space of \emph{isomorphisms} between elements $x,y : C$ in a precategory is defined by
\[
x \cong y \coloneqq \Sigma_{f : \Hom(x,y)} \isIso(f) \qquad \text{where} \qquad \isIso(f) \coloneqq  \left( \Sigma_{g : \Hom(y,x)}g \circ f \sim \id_x \right) \times \left( \Sigma_{h : \Hom(y,x)} f \circ h \sim \id_y\right) \, {.}
\] 
The unexpected complexity in the definition of $\isIso(f)$ is because in the higher setting we must use a homotopy coherent model of the \emph{free living isomorphism} $\iso$ in $\sSet$.

  \begin{definition}\label{defn:rezk}

    A precategory is a \textbf{category} if for all $x,y : C$, the natural map is an equivalence: 
    \[ \begin{tikzcd} C \arrow[dr, "{(x \sim y)_{x,y :C}}"'] \arrow[rr, "\pathind(\refl\mapsto \id)"] & & C^\iso \arrow[dl, "{(x \cong y)_{x,y : C}}"] \\ & C \times C\, {.} \end{tikzcd}\] 
  \end{definition}

  \begin{definition}\label{defn:discrete} 
    A simplicial space $C$ is a \textbf{groupoid} if for all $x,y : C$, the natural map is an equivalence: 
    \[ \begin{tikzcd} C \arrow[dr, "{(x \sim y)_{x,y :C}}"'] \arrow[rr, "\pathind(\refl\mapsto \id)"] & & C^{\Delta^1} \arrow[dl, "{(\Hom(x, y))_{x,y : C}}"] \\ & C \times C\, {.} \end{tikzcd}\] 
  \end{definition}

One can show that a groupoid is a category in which all arrows are isomorphisms, meaning the restriction $C^\iso \fto C^{\Delta^1}$ defined by the family $(\isIso(f))_{f : C^{\Delta^1}}$ is an equivalence.

\section{Arrow induction.}\label{sec:arrow}

Arrow induction extends the Yoneda lemma from traditional category theory, which describes the mapping out universal property of a representable functor. To use more familiar language, we refer to ``categories'' here but all the results about categories in this section actually apply in the more general setting of precategories, satisfying only \cref{defn:segal}.

There are various formulations of the Yoneda lemma for families of categories first described by Street \cite{Street1980fb}, then adapted to the higher categorical setting by \cite{RiehlVerity2017fy, RiehlVerity2022}, and to the setting of synthetic categories by \cite{RiehlShulman2017,BuchholtzWeinberger2023sf,Weinberger2024ts}. Here we present the variant that is easiest to explain, rather than the most general, developing the Yoneda lemma for families of $\infty$-groupoids that vary covariantly functorially in arrows in the base $\infty$-category. Other variants consider families of $\infty$-categories which might mix covariant and contravariant dependencies.

\subsection{Arrow induction for covariant families.}

Fixing a category $C$ and element $c : C$, the covariant representable functor is encoded by the family of simplicial spaces $(\Hom(c,x))_{x : C}$. This family has a special property of being \emph{covariant}, meaning any arrow $f \colon \Hom(x,y)$ in the base simplicial space lifts uniquely to an arrow in the total space with specified source.

\begin{definition}\label{defn:covariant}
  A family of simplicial spaces $(E_x)_{x : B}$ is \textbf{covariant} if for all arrows $f : \Hom(x,y)$ in $B$ the function $\ev_0 \colon \Pi_{t : \Delta^1} E_{f(t)} \to E_x$ is an equivalence.
\end{definition}
Here $\Pi_{t : \Delta^1} E_{f(t)}$ is the simplicial space of arrows in the family living over the arrow $f$ in the base. \cref{defn:covariant} asserts that such arrows are uniquely determined by their source, or equivalently that the square
\begin{equation}\label{eq:covariant} \begin{tikzcd} E^{\Delta^1} \arrow[r, two heads, "\ev_0"] \arrow[d, two heads, "{(\Pi_{t : \Delta^1} E_{f(t)})_{f : B^{\Delta^1}}}"'] \arrow[dr, phantom, "\lrcorner" very near start]& E \arrow[d, two heads, "{(E_x)_{x : B}}"] \\ B^{\Delta^1} \arrow[r, two heads, "\ev_0"'] & B \end{tikzcd} \qquad \leftrightsquigarrow \qquad \begin{tikzcd} \Pi_{t : \Delta^1} E_{f(t)} \arrow[r, uwearrow, "\ev_0"] & E_x \end{tikzcd} \text{for~all~} f: \Hom(x,y)~\text{in~}B \end{equation} 
is a pullback in the higher categorical sense. By a contractibility argument:

\begin{lemma}\label{lem:groupoid-fibers}
If $(E_x)_{x :B}$ is a covariant family of simplicial spaces, then the fibers $E_x$ are groupoids.\footnote{The original reference \cite[8.18]{RiehlShulman2017} assumes that the base $B$ is a precategory, but this hypothesis is unnecessary, as has been verified in an accompanying library where results about synthetic categories have been formalized in the computer proof assistant \textsc{Rzk} \cite{sHoTT}.} 
\end{lemma}

We now explain the terminology ``covariant.''
For each arrow $f : \Hom(x,y)$ in the base of a covariant family, there is a map between the fibers $f_* \colon E_x \to E_y$ defined by composing the inverse of the equivalence of \cref{defn:covariant} with the other projection $\ev_1 \colon \Pi_{t : \2} E_{f(t)} \to E_y$. As with \cref{defn:segal}, the advantage of formulating \cref{defn:covariant} as we did is that it follows that these actions are covariantly functorial when the base is a category:

\begin{lemma}
  When $B$ is a category and $(E_x)_{x : B}$ is a covariant family of simplicial spaces, the action of arrows in the base on the fibers is covariantly functorial: we have $\id_x(u) \sim u$ for $u : E_x$ and $(g \circ f)_* u \sim g_* (f_* u)$ for $u : E_x$, $f : \Hom(x,y)$, and $g : \Hom(y,z)$.
\end{lemma}

Just as functoriality comes for free for synthetic categories, naturality comes for free for maps between covariant families: 

\begin{lemma}\label{lem:nat-for-free}
  Given covariant families $(E_x)_{x :B}$ and $(F_x)_{x : B}$, any fiberwise map $\phi : \Pi_{x : B} (E_x \to F_x)$ as below-left is automatically natural: for any $f : \Hom(x,y)$, the diagram below-right commutes up to homotopy.
  \[ 
  \begin{tikzcd} E \arrow[rr, "\phi"] \arrow[dr, two heads, "{(E_x)_{x : B}}"'] &[-2em] & [-2em] F \arrow[dl, two heads, "{(F_x)_{x : B}}"] & & & E_x \arrow[d, "f_*"'] \arrow[r, "\phi_x"] & F_x \arrow[d, "f_*"] \\ & B & && & E_y \arrow[r, "\phi_y"'] & F_y\, {.} \end{tikzcd}\] 
\end{lemma}

\begin{lemma} A simplicial space $C$ is a precategory if and only if for each $c : C$, the family $(\Hom(c,x))_{x : C}$ is covariant, in which case the action on fibers $f_* \colon \Hom(c,x) \to \Hom(c,y)$ coincides with composition $f \circ - \colon \Hom(c,x) \to \Hom(c,y)$ with $f : \Hom(x,y)$.\footnote{The original published proof of this result \cite[8.13]{RiehlShulman2017} had an error of circular reasoning that was uncovered in the process of formalizing the proof in \texttt{Rzk} \cite{KudasovRiehlWeinberger2024}.}
\end{lemma}

We can now state the arrow induction principle.

\begin{proposition}[arrow induction]\label{prop:arrow-ind}
  Fix a category $C$, an element $c : C$, and consider the covariant family $(\Hom(c,x))_{x : C}$ whose total space is the coslice category $\coslice{c}{C}$ of arrows under $c$. 
  To define a section of a covariant family over the coslice category $\coslice{c}{C}$,  it suffices to define the image of the identity arrow $\id_c$ in the fiber $F_{\id_c}$.
  \[ \begin{tikzcd}[column sep=huge] \Delta^0 \arrow[d, "\id_c"'] \arrow[r, "d"]  & F \arrow[d, two heads, "{(F_f)_{f : \coslice{c}{C}}}"] \\ \coslice{c}{C} \arrow[r, equals] \arrow[ur, dashed, "{\arrowind(\id_c \mapsto d)}" description] & \coslice{c}{C}\, {.} \end{tikzcd} \] 
\end{proposition}

The semantic justification for \cref{prop:arrow-ind} is similar to the justification for \cref{prop:basic-path-ind}. The functor $\id_c \colon \Delta^0 \to \coslice{c}{C}$ is \emph{initial}, which in this instance just means that $\id_c : \coslice{c}{C}$ is an initial element. Initial functors $i \colon A \to X$ have a left lifting property against covariant families $(E_b)_{b : B}$, which may be internalized to the statement that the induced map on mapping spaces is a trivial fibration 
\[
    \begin{tikzcd} A \arrow[d, "i"'] \arrow[r, "e"] & E \arrow[d, two heads, "{\rho = (E_b)_{b : B}}"] & &\arrow[d, phantom, "\leftrightsquigarrow"] & \Map(X, E) \arrow[d, utfibarrow, "\rho \circ - \circ i"] \\ X \arrow[r, "b"'] \arrow[ur, dashed] \arrow[ur, dashed] & B &&  ~ & \Map(A,E) \times_{\Map(A,B)} \Map(X,B)\, {.}  \arrow[u, bend left, dashed]
      \end{tikzcd}        
\]
When expressed internally to the language of synthetic categories, arrow-induction provides a section to the map
  \[
\begin{tikzcd} \Delta^0 \arrow[d, "\id_c"'] \arrow[r, "d"]  & F \arrow[d, two heads, "{\rho = (F_f)_{f : \coslice{c}{C}}}"] \\ \coslice{c}{C} \arrow[r, equals] \arrow[ur, dashed] & \coslice{c}{C} \end{tikzcd}  
 \qquad
    \begin{tikzcd} 
      \Map_{\coslice{c}{C}}(\coslice{c}{C},F) \arrow[d, utfibarrow, "- \circ \id_c"] \arrow[r] \arrow[dr, phantom, "\lrcorner" very near start] & \Map(\coslice{c}{C},F) \arrow[d, utfibarrow, "\rho \circ - \circ \id_c"] \\     
       \Map(\Delta^0, F_{\id_c}) \arrow[r]   \arrow[u, dashed, bend left, "\arrowind"] \arrow[d, two heads] \arrow[dr, phantom, "\lrcorner" very near start] & \Map(\Delta^0,F) \times_{\Map(\Delta^0,\coslice{c}{C})} \Map(\coslice{c}{C}, \coslice{c}{C})  \arrow[u, dashed, bend left] \arrow[d, two heads, "\pi"] \\
       \Delta^0 \arrow[r, "\id"'] & \Map(\coslice{c}{C}, \coslice{c}{C})\, {.}
    \end{tikzcd}
\]
Arrow induction defines an equivalence between the sections of a covariant family $(F_f)_{f : \coslice{c}{C}}$ and the fiber $F_{\id_c}$:
\[ \begin{tikzcd} \Map_{\coslice{c}{C}}(\coslice{c}{C},F) \arrow[r, utfibarrow, "\ev_{\id_c}"] & F_{\id_c} \, {.}\arrow[l, bend left, uwedashed, "\arrowind"] \end{tikzcd}\] 

\subsection{The Yoneda Lemma.}

In \cite{RiehlShulman2017}, \cref{prop:arrow-ind} is referred to as the \emph{dependent Yoneda lemma}. The more familiar Yoneda lemma can be derived as a corollary. Consider a category $C$, an element $c :C$, and a covariant family $(F_x)_{x :C}$. By \cref{lem:nat-for-free}, $\Pi_{x :C} (\Hom(c,x) \to F_x)$ defines the simplicial space of natural transformations from the first family to the second. There is a map from the simplicial space of natural transformations to $F_c$ that sends $\phi$ to $\phi_c(\id_c)$ that we refer to as $\evid$.

\begin{proposition}[Yoneda lemma]\label{prop:yoneda} For a category $C$, an element $c : C$, and a covariant family $(F_x)_{x :C}$, the natural map is an equivalence:
  \[ \begin{tikzcd} \left( \Pi_{x :C} (\Hom(c,x) \to F_x)  \right) \arrow[r, "\evid"]& F_c\, {.}\end{tikzcd} \] 
\end{proposition}

In the internal language for synthetic categories, \cref{prop:yoneda} is arguably both simpler to state and prove than the Yoneda lemma for strict 1-categories \cite{Riehl2023u}. In our semantic setting, the argument is as follows:

\begin{proof}
In the following diagram, note that by the universal property of the pullback $\Map_{\coslice{c}{C}}(\coslice{c}{C},P) \cong \Map_C(\coslice{c}{C},F)$ and $P_{\id_c} \cong F_c$. Thus, the top left-hand map below-right encodes the desired equivalence:
  \[
    \begin{tikzcd}[column sep=large]
      \Delta^0 \arrow[rr, bend left, "e"]\arrow[d, tcofarrow, "\id_c"'] \arrow[r, "d"] & P \arrow[d, two heads, "\chi"] \arrow[r, two heads] \arrow[dr, phantom, "\lrcorner" very near start] & F \arrow[d, two heads, "{(F_x)_{x : C}}"] & \arrow[d, phantom, "\leftrightsquigarrow"]   \\  
      \coslice{c}{C} \arrow[r, equals] \arrow[ur, dashed] & \coslice{c}{C} \arrow[r, "{(\Hom(c,x))_{x : C}}"', two heads] & C & ~   \end{tikzcd}   \quad
    \begin{tikzcd} 
      \Map_{\coslice{c}{C}}(\coslice{c}{C},P) \arrow[d, utfibarrow, "- \circ \id_c"] \arrow[r] \arrow[dr, phantom, "\lrcorner" very near start] & \Map(\coslice{c}{C},P) \arrow[d, utfibarrow, "\chi \circ - \circ \id_c"] \\     
       \Map(\Delta^0, P_{\id_c}) \arrow[r]   \arrow[u, dashed, bend left, "\arrowind"] \arrow[d, two heads] \arrow[dr, phantom, "\lrcorner" very near start] & \Map(\Delta^0,P) \times_{\Map(\Delta^0,\coslice{c}{C})} \Map(\coslice{c}{C}, \coslice{c}{C})  \arrow[u, dashed, bend left] \arrow[d, two heads, "\pi"] \\
       \Delta^0 \arrow[r, "\id"'] & \Map(\coslice{c}{C}, \coslice{c}{C})\, {.}
    \end{tikzcd} 
\]
\end{proof}

\section{Directed univalent universes.}\label{sec:dua}

Simplicial type theory is an efficient formal system to state and prove theorems about categories \cite{BardomianoMartinez2025lc}. But the original formulation of \cite{RiehlShulman2017} does not provide any techniques for building examples of synthetic categories. As we saw in \S\ref{ssec:ua-applications}, interesting examples of synthetic spaces can be built from the universal family of spaces. In this section, we will sketch the construction of an analogous universal covariant family of simplicial spaces. As a consequence of a \emph{directed} analog of Voevdosky's univalence principle, the base of the universal covariant family is a category. By \cref{lem:groupoid-fibers}, its fibers are groupoids, aka ``spaces'' in our present terminology, so this defines the category of spaces, which can be used to construct other interesting examples.

In \S\ref{ssec:dua}, we adapt the approach of \cite{Shulman2019} to construct a directed univalent universe of simplicial spaces---or more generally, in a strict 1-categorical presentation of the category of simplicial objects in any $\infty$-topos---previewing results from forthcoming joint work with Cavallo and Sattler \cite{CavalloRiehlSattler}.  There are various other approaches to constructing a directed univalent universe of spaces or categories \cite{NicholsBarrer2007, Lurie2009, KazhdanVarshavsky2014, Cisinski2019, WeaverLicata2020, CisinskiNguyen2022,GratzerWeinbergerBuchholtz:2024du}. In \S\ref{ssec:dua-applications}, we very briefly sketch applications; see \cite{CisinskiNguyen2022, CisinskiCnossenNguyenWalde, GratzerWeinbergerBuchholtz:2024du} for more. 

\subsection{Constructing directed univalent universes.}\label{ssec:dua}

As discussed in \S\ref{ssec:uu}, 
our objective is to define a universal covariant family of simplicial spaces so that $\Gamma$-indexed generalized elements define covariant families of small simplicial spaces in context $\Gamma$. We work in the strict 1-category of simplicial objects in a type theoretic model topos, for which the 1-category $\sSet^{\DDelta^\op}$ introduced in \S\ref{sec:cats} is the basic example; we continue to refer to generic (fibrant) objects as simplicial spaces. In that setting, our covariant families are more commonly called \emph{left fibrations}, these being maps $\rho \colon E \fto B$ that define fibrations in the ambient model structure with the property that the induced map to the strict 1-categorical pullback in the commutative square \cref{eq:covariant} is a trivial fibration. We show that suitably structured left fibrations define a notion of fibred structure that is small groupoid valued and locally representable, and thus admits a universe $\varrho \colon \cS_* \fto \cS$ that is univalent, satisfying the condition of \cref{defn:univalence}. The covariant universe $\cS$ has a classifying property analogous to \cref{eq:univ-pullback} that implies that any small covariant family $(A_\gamma)_{\gamma : \Gamma}$ has an ``straightening'' into a functor $A \colon \Gamma\to \cS$ whose pullback recovers the  ``unstraightened'' covariant family of small groupoids 
\[ 
\begin{tikzcd} A \arrow[d, two heads, "{(A_\gamma)_{\gamma :\Gamma}}"'] \arrow[dr, phantom, "\lrcorner" very near start] \arrow[r] & \cS_\bullet \arrow[d, two heads, "\varrho"] \\ \Gamma \arrow[r, "A"'] & \cS\, {.}
\end{tikzcd}
\] 
By univalence of $\varrho \colon \cS_\bullet \fto \cS$, equivalent covariant families over $\Gamma$ are classified by homotopic functors into $\cS$.

The directed univalence axiom characterizes the arrows $(\Hom(A,B))_{A,B : \cS}$ in $\cS$ in terms of another universal family. There is a family of simplicial spaces $(A \to B)_{A, B :\cS}$ defined by using the internal hom in the category of simplicial spaces over $\cS \times \cS$ \cite[6.2.2]{Stenzel2019phd}. By construction a generalized element 
\[
  \begin{tikzcd}  
    {(A \to B)} \arrow[r] \arrow[dr, phantom, "\lrcorner" very near start] \arrow[d, two heads] & \Fun(\cS_\bullet) \arrow[d, two heads, "{(A \to B)_{A,B : \cS}}"] \\ \Gamma \arrow[u, bend left, dashed, "{e}"] \arrow[r, "{(A,B)}"']& \cS \times \cS
  \end{tikzcd} 
  \]
encodes a function $f : A \to B$ between the covariant families over $\Gamma$ classified by the maps $A,B \colon \Gamma \to \cS$.

\begin{definition}\label{defn:directed-univalence}
  The fibration $\varphi \colon \cS_* \fto \cS$ is \textbf{directed univalent} if the comparison map from the space of arrows to the space of functions defined by arrow induction\footnote{Here the arrow induction principle of \cref{prop:arrow-ind} should be interpreted in context $\cS$, holding the variable $A: \cS$ fixed.} is an equivalence:
\[
  \begin{tikzcd}[column sep=huge] 
    \cS\arrow[d, "\id"', tail] 
    \arrow[r, "\id"] & 
    {\Fun(\cS_\bullet)} \arrow[d, two heads,"{(A \to  B)_{A,B: \cS}}"] \\ 
    \cS^{\Delta^1} \arrow[ur, dashed, "{\arrfun}" description] \arrow[r, two heads, "{(\Hom(A,B))_{A, B: \cS}}"'] & \cS \times \cS
  \end{tikzcd} \qquad \arrfun \coloneqq \arrowind(\id \mapsto \id)\, {.}
  \]
\end{definition}

See \cite{CavalloRiehlSattler} for more details and a proof that the universal covariant family $\varrho \colon \cS_\bullet \fto \cS$ is directed univalent. There we also establish a ``higher'' version of directed univalence defining an equivalence between the simplicial space  of $n$ composable arrows in $\cS$ and the simplicial space of $n$ composable functions between covariant families over $\cS$. We also show that the functorial action in the universal covariant family is given by function composition.

\subsection{Applications of directed univalence.}\label{ssec:dua-applications}

The higher version of directed univalence can be used to provide the first example of a synthetic category.

\begin{corollary} The base of the universal covariant family defines a category.
\end{corollary}

As noted above, by \cref{lem:groupoid-fibers}, the points in $\cS$ are groupoids, aka spaces, so we refer to the base of the universal covariant family as the \emph{category of spaces}.
The total space $\cS_\bullet$ of the universal covariant family is a category as well, namely the category of pointed spaces $(A,a)$, explaining our notation.
\[ \begin{tikzcd} A \arrow[d, two heads] \arrow[r] \arrow[dr, phantom, "\lrcorner" very near start] & \cS_\bullet \arrow[d, two heads, "\varrho"] \\ \Delta^0 \arrow[r, "A"'] \arrow[u, dashed, bend left, "a"] & \cS\, {.}\end{tikzcd}\]

A useful consequence of univalence, mentioned in \S\ref{ssec:ua-applications}, is the structure identity principle, characterizing path spaces of various classifying spaces built from the universe. Similarly, directed univalence implies a directed analog, which we might call the \emph{structure homomorphism principle}, identifying the arrows in various classifying categories built from the universe as homomorphisms of the appropriate structures \cite{WeaverLicata2020,GratzerWeinbergerBuchholtz:2024du}. For example, arrows in the category $\cS_\bullet$ of pointed spaces can be identified with pointed functions. 

\begin{example}
 We can define the simplicial space of \emph{magmas} in $\cS$ by 
 \[ \type{Magma} \coloneqq \Sigma_{A : \cS} A \times A \to A .\] 
 The elements are spaces $A$ equipped with a binary operation $\mu \colon A \times A \to A$. By directed univalence, $\type{Magma}$ is a category whose arrows are homomorphisms. For other examples, see \cite[\S6]{GratzerWeinbergerBuchholtz:2024du}.
 \end{example}

\section{Vistas.}\label{sec:vistas}

To conclude we highlight ongoing work involving related ideas.

\subsection{Computational homotopy type theory.}\label{ssec:hott}

While homotopy type theory is now almost two decades old, the precise rules of this formal system are still under development. One motivation for further innovation is to restore some of the metatheoretical properties of intensional dependent type theory, which as a foundation for constructive mathematics satisfies canonicity and normalization. In such a system a proof establishing the existence of an integer $n$ with a certain property---such as famously Brunerie's proof that $\pi_4S^3$ is $\mathbb{Z}/n\mathbb{Z}$ for some $n : \mathbb{Z}$---can be compiled out to calculate the value of $n$ \cite{Brunerie2019, LjungstromMortberg2023}.  This motivated various \emph{cubical} variants of homotopy type theory surveyed in \cite{CavalloMortbergSwan2020}, 
though there are still some open questions concerning the precise homotopy theoretic semantics of such systems \cite{AwodeyCavalloCoquandRiehlSattler2024, CavalloSattler2025}.

A related objective is to convert the equivalences characterizing identity types---such as those given by the structure-identity principle---into definitional equalities. Such a system, called \emph{Higher Observational Type Theory}, is currently under development by Altenkirch, Kaposi, and Shulman, extending prior joint work on internal parametricity \cite{AltenkirchChamounKaposiShulman}, and has an accompanying experimental proof assistant \textsc{Narya}.

\subsection{Multimodal simplicial type theory.}\label{ssec:mmtt}

Theorems about spaces that are not equivalence-invariant cannot be correctly stated or proven in homotopy type theory for the reasons discussed in \S\ref{ssec:ua-applications}. Similarly, constructions involving categories that are not functorial or natural cannot be correctly stated or prove in simplicial type theory by \cref{lem:fun-for-free,lem:nat-for-free}. Thus, further axioms or extensions of this formal system are needed to redevelop the full analytic theory of categories in a synthetic setting.

Work in progress towards this goal is being pursued by Cisinski, Cnossen, Nguyen, and Walde \cite{CisinskiCnossenNguyenWalde} and Gratzer, Weinberger, and Buchholtz \cite{GratzerWeinbergerBuchholtz:2024du, GratzerWeinbergerBuchholtz:2025ye}.
The latter project extends the simplicial type theory with various \emph{modalities}, which, for instance, might distinguish between contexts that are categorical or groupoidal and limit the use of variables in certain circumstances.

\subsection{Synthetic mathematics.}\label{ssec:synthetic}

Synthetic mathematics, in the tradition of Euclid, is experiencing a revival, perhaps in relation to recent advances in computer proof assistants, which can help mathematicians trained in one formal framework learn to reason in another \cite{Shulman2024sf}.
A new proposed synthetic approach to algebraic geometry by Blechschmidt and Cherubini, Coquand, and Hutzler provides axioms internally to the Zariski topos as a higher topos, so as to have direct access to higher homotopy types to reason about cohomology \cite{Blechschmidt2017, CherubiniCoquandHutzler2024}. There are similar new proposed synthetic approaches to condensed mathematics intended to provide internal languages for the 1-topos of light condensed sets or the higher topos of light condensed anima developed in forthcoming work of Barton and Commelin and by Cherubini, Coquand, Geerligs, and Moeneclaey \cite{CherubiniCoquandGeerligsMoeneclaey2025}.

\section*{Acknowledgments.}

Part of this text was originally prepared as lecture notes for a talk given at ``A panorama of homotopy theory: a conference in honor of Mike Hopkins.'' Expository suggestions came from conversations with Ulrik Buchholtz, Evan Cavallo, Daniel Gratzer, Egbert Rijke, Mike Shulman, and Jonathan Weinberger, while  Katsunori Takahashi caught a typo. 
 The author is grateful for support from the US National Science Foundation via the grant DMS-2507077 and the US Air Force Office of Scientific Research under award number FA9550-21-1-0009.

\bibliographystyle{siamplain}
\bibliography{synthetic.bbl}
\end{document}